\documentclass[11pt]{article}
\usepackage{latexsym}
\usepackage[all]{xy}
\usepackage{amsmath}
\usepackage{amssymb}
\usepackage{amsfonts}
\usepackage{stmaryrd}

\newtheorem{thm}{Theorem}[subsection]
\newtheorem{prop}[thm]{Proposition}
\newtheorem{lem}[thm]{Lemma}
\newtheorem{df}[thm]{Definition}
\newtheorem{cor}[thm]{Corollary}
\newtheorem{conj}[thm]{Conjecture}

\newtheorem{rmk}[thm]{Remark}

\newtheorem{Q}[thm]{Question}

\begin{document}

\title{\textbf{Segal topoi and stacks over Segal categories}}
\bigskip
\bigskip

\author{\bigskip\\
Bertrand To\"en \\
\small{Laboratoire J. A. Dieudonn\'{e}}\\
\small{UMR CNRS 6621} \\
\small{Universit\'{e} de Nice Sophia-Antipolis}\\
\small{France}\\
\bigskip \and
\bigskip \\
Gabriele Vezzosi \\
\small{Dipartimento di Matematica}\\
\small{Universit\`a di Bologna}\\
\small{Italy}\\
\bigskip}

\date{}

\maketitle

\begin{abstract}
In \cite{partI} we began the study of higher sheaf theory (i.e. stacks theory)
on higher categories endowed with a suitable notion of topology:
precisely, we defined the notions of $S$\textit{-site} and of
\textit{model site}, and the associated \textit{categories of stacks}
on them. This led us to a notion of \textit{model topos}. In this paper we treat the analogous
theory starting from (1-)\textit{Segal categories} in place of $S$-categories and model categories.

We introduce notions
of \textit{Segal topologies}, \textit{Segal sites} and stacks over them, giving rise
to a definition of \textit{Segal topos}. We compare the notions of Segal topoi and
of model topoi, showing that the two theories are equivalent in some sense. However,
the existence of a nice Segal category of morphisms between Segal
categories allows us to improve the treatment of topoi in this context. In particular
we construct the \textit{2-Segal category of Segal topoi and
geometric morphisms}, and we provide a Giraud-like statement charaterizing Segal topoi
among Segal categories.

As example of applications, we show how to reconstruct a topological space from
the Segal topos of \textit{locally constant stacks} on it, thus extending the
main theorem of \cite{to3} to the case of un-based spaces. We also give some
hints of how to define \textit{homotopy types of Segal sites}: this approach gives a new
point of view and some improvements on the \'etale homotopy theory of schemes, and more
generally on the theory of homotopy types of
Grothendieck sites as defined by Artin and Mazur.
\end{abstract}


\medskip

\textsf{MSC-class:} $14$A$20$; $18$F$10$; $55$U$40$; $18$D$20$; $18$D$05$.

\tableofcontents

\bigskip

\begin{section}{Introduction}

In \cite{partI} we have developed a homotopy version of sheaf theory, in which not only sheaves of sets are
replaced by certain simplicial presheaves (called \textit{stacks}), but in which, more crucially, the base Grothendieck sites are
replaced either by \textit{$S$-sites}, i.e. simplicially enriched categories endowed with a suitable notion
of topology, or by \textit{model sites}, i.e. model categories endowed with a suitable notion of topology.
These constructions have led us naturally to the concept of \textit{model topos}, originally due to
C. Rezk (\cite{rez2}), for which one can prove some generalizations of the basic results in topos theory,
as for example the correspondence between topologies and certain localizations of categories of presheaves
(see \cite[Theorem 3.8.3]{partI}). However, it seems very uneasy to further develope the theory of
model topoi purely in terms of model categories, and one reason for this
is the missing theory of internal Hom-objects, or equivalently the non-existence of
reasonable model categories of Quillen functors between two model categories. For example,
it seems very difficult to construct an analog for model topoi of the $2$-category of topoi and geometric
morphisms between them. The main purpose of this work is to solve this problem by introducing similar constructions
in the context of \textit{Segal categories} of \cite{sh,p}. \\

\textbf{Segal categories, Segal sites and stacks.} Segal categories are
weak form of $S$-categories, and behave very much the same way. In fact the theory of
Segal categories and of $S$-categories are equivalent in some sense (see \S 4.1), and
it is not a bad idea to think of Segal categories as $S$-categories, at least as
a first approximation. Thanks to the fundational work of C. Simpson, A. Hirschowitz and R. Pellissier (see
e.g. \cite{sh,p,s3}), many (if not all) standard categorical notions are available
for Segal categories, as for example categories of functors, adjunctions, limits and
colimits, Yoneda lemma \dots (see \S 3.1, 3.2, 4.1, 4.2). These constructions will allow us to
follow the main line of topos theory in the new context of Segal categories.

We define a \textit{Segal topology} on a Segal category $T$ to be a Grothendieck topology on its
homotopy category $\mathrm{Ho}(T)$ (see Definition \ref{dtop}). A Segal category endowed with
a topology will be called a \textit{Segal site}. For a Segal site $(T,\tau)$, one defines
a Segal category of \textit{pre-stacks} $\hat{T}$, as well as a sub-Segal category
of \textit{stacks} $T^{\sim,\tau}$ (see Definition \ref{dstacks} (3)). \\

\textbf{Segal topoi and Segal category of stacks.}
A \textit{Segal topos} is defined to be a Segal category equivalent to some
exact localization of a Segal category of pre-stacks $\hat{T}$ (see
Definition \ref{dsegtop}).
If $(T,\tau)$ is a Segal site, the sub-Segal category
$T^{\sim,\tau}$ is a left exact localization of $\hat{T}$, and therefore is
a Segal topos. Our first result, Theorem \ref{cmp} states that
a large class of Segal topoi (the ones being \textit{$t$-complete}) are obtained this way.
In other words, the map $\tau \mapsto T^{\sim,\tau}$ gives
a one-to-one correspondence between topologies on $T$ and
left exact localizations of $\hat{T}$ which are furthermore $t$-complete.
The latter result is a generalization to the case of Segal topoi of the well known
fact that topologies on a category are in bijection with left exact localizations
of its category of presheaves. It also justifies our definition of a topology
on a Segal category. We also provide a Giraud's style statement
charaterizing Segal topoi among Segal categories (Conjecture \ref{giraud}). \\

\textbf{Geometric morphisms and the $2$-Segal category of Segal topoi.}
A crucial property of Segal categories is the existence of a
reasonable Segal category of morphisms between two Segal categories.
Using this construction, as well as the notion of adjunction and limits,
we define for two Segal topoi $A$ and $B$
a Segal category of \textit{geometric morphisms} $\mathbb{R}\underline{Hom}^{geom}(A,B)$
(see Definition \ref{dgeom}). These Segal categories of geometric morphisms
assemble together into a \textit{$2$-Segal category} $\mathsf{SeT}$, of Segal
topoi, which is our Segal analog of the $2$-category of topoi and geometric morphisms. \\

\textbf{Segal topoi and model topoi.} Let $(T,\tau)$ be an $S$-site (i.e.
a Segal site where $T$ is furthermore an $S$-category). On one hand
there is the model category of stacks $\mathsf{SPr}_{\tau}(T)$, constructed in
\cite[Theorem 3.4.1]{partI}, and on the other there is
the Segal topos $T^{\sim,\tau}$ of stacks over $(T,\tau)$ defined
in Definition \ref{dstacks}. Using the strictification theorem
proved in \cite{sh}, we show that
there exists a natural equivalence of Segal categories (see Corollary \ref{ccomp})
$$L(\mathsf{SPr}_{\tau}(T))\simeq T^{\sim,\tau},$$
where $LM$ is the simplicial localization of a model category $M$ as defined
in \cite{dk1}. This shows that the construction $M\mapsto LM$ gives
a way to pass from model topoi to Segal topoi. Therefore, if $M$ and $N$ are two model
topoi, one can consider $\mathbb{R}\underline{Hom}^{geom}(LM,LN)$, the Segal
category of geometric morphisms from $LM$ to $LN$. This construction gives a solution
to our original problem of defining a ``category'' of geometric morphisms between two model topoi. \\

\textbf{An application: Galois interpretation of homotopy theory}.
As an example of application of our notion of Segal categories of geometric morphisms
we provide a \textit{Galois interpretation of homotopy theory}, extending
the well known relations between fundamental groupoids and categories of
locally constant sheaves on a space. For a $CW$ complex $X$, one consider the Segal category $\mathsf{Loc}(X)$ of locally constant stacks on $X$, which is
easily seen to be a Segal topos. For two $CW$ complexes $X$ and $Y$ we prove that
$\mathbb{R}\underline{Hom}^{geom}(\mathsf{Loc}(X),\mathsf{Loc}(Y))$ is in fact a Segal groupoid whose
geometric realization is equivalent to $\mathbb{R}\underline{Hom}(X,Y)$
(see Theorem \ref{tgal}). In particular, taking $X=*$ we get that
$Y$ is weakly equivalent to the geometric realization of the Segal groupoid
$\mathbb{R}\underline{Hom}^{geom}(\mathsf{Top},\mathsf{Loc}(Y))$ (here
$\mathsf{Top}:=L\mathsf{SSet}$ is the Segal category of simplicial sets). This last statement generalizes
to the whole homotopy type the fact that
the fundamental groupoid of a space is equivalent to the groupoid of fiber functors
on its category of locally constant sheaves (see \cite[Exp. V]{sga1}). \\

\textbf{Homotopy types of Segal topoi.} Based on our
Galois interpretation of homotopy theory we
define for any Segal topos $T$ a morphism of Segal categories
$H_{T} : \mathsf{Top} \longrightarrow \mathsf{Top}$, which has to be thought as some
kind of homotopy type of $T$; we call $H_{T}$ the
\textit{homotopy shape of $T$} (see Definition \ref{dhom}). The morphism $H_{T}$
is in fact pro-corepresentable by a pro-homotopy type
$X$. This pro-homotopy type (which is a pro-object in $\mathsf{Top}$ rather than a
pro-object in $\mathrm{Ho}(\mathsf{Top})$) is called the \textit{homotopy type of the Segal topos $T$}.
When $T$ is the Segal topos of stacks over a Grothendieck site, this approach gives a new point of
view on homotopy types of sites as defined by Artin and Mazur (\cite{am}). This also allows us
to define the \'etale pro-homotopy type of the sphere spectrum which seems to us
an interesting example to consider (see \ref{quest2}).  \\

\bigskip
\bigskip

\noindent \textbf{Related works.} In \cite{s3} C. Simpson
investigates \textit{Segal pre-topoi} and the question of the existence of a theory of Segal topoi is at
least implicit in the text, if not clearly stated. In some sense this work, and more
precisely our conjecture \ref{giraud} gives a possible answer to his question.

As we have already mentioned,
in his unpublished manuscript \cite{rez2} C. Rezk has introduced a notion of \textit{homotopy
topos}, which is a model category version of our definition of Segal topos. More recently,
J. Lurie has investigated a notion of $\infty$-topos, equivalent to our notion
of Segal topos, and for which he proved a Giraud's theorem (therefore
our conjecture \ref{giraud} seems to be a theorem now). He has also related the notion
of $\infty$-topos with a different notion of stacks (see \cite{lu}).

Let us also mention that
A. Joyal (see \cite{jo2}) has developed a theory of \textit{quasi-categories}, which is expected to be
equivalent to the theory of $S$-categories and of Segal categories, and for which he has
defined a notion of \textit{quasi-topos} very similar to
our notion of Segal topos. The two approaches are expected to be
equivalent.

The work of D-C. Cisinski \cite{cis} seems to be closely
related to a notion of \textit{hypertopology} we discuss in Remark \ref{hyper}.

In his letter to L. Breen (see \cite{gr}), A. Grothendieck
proposes some Galois interpretation of homotopy types of sites. We like
to consider our Theorem \ref{tgal} as a topological version of this
(see also the Remark after Corollary 3.2 in \cite{to3}). \\

\bigskip

\noindent \textbf{Acknowledgments.} We wish to thank the MSRI in Berkeley
for invitation to the Program ``Stacks, Intersection Theory and Non Abelian Hodge Theory'' (January-May 2002).\\
In particular we wish to acknowledge many fruitful discussions with the
participants to the above program. We especially thank A.
Joyal, M. Kontsevich, J. Lurie, and C. Rezk for useful
exchanges. We are particularly indebted to C. Simpson for suggestions, discussions
on our project and for all his deep work on higher categories and higher stacks.\\

\noindent \textbf{Notations and conventions.}
We will use the word \textit{universe} in the sense of \cite[Exp. I, Appendice]{sga4}.
We will fix three universes $\mathbb{U} \in \mathbb{V} \in \mathbb{W}$, and assume
that $\mathbb{N} \in \mathbb{U}$. The category of sets (resp. simplicial sets, resp. \dots)
belonging to a universe $\mathbb{U}$ will be denoted by $\mathsf{Set}_{\mathbb{U}}$ (resp. $\mathsf{SSet}_{\mathbb{U}}$, resp. \dots).
The objects of $\mathsf{Set}_{\mathbb{U}}$ (resp. $\mathsf{SSet}_{\mathbb{U}}$, resp. \dots) will be called
$\mathbb{U}$-sets (resp. $\mathbb{U}$-simplicial sets, resp. \dots). We will use the expression
\textit{$\mathbb{U}$-small set}
(resp. \textit{$\mathbb{U}$-small simplicial set}, resp. \dots) to mean \textit{a set isomorphic
to a set in $\mathbb{U}$}
(resp. \textit{a simplicial set isomorphic to a simplicial set in $\mathbb{U}$}, resp. \dots).

Our references for model categories are \cite{ho} and \cite{hi}.
By definition, our model categories will always be \textit{closed} model categories, will have all \textit{small}
limits and colimits and the functorial factorization property.
The word \textit{equivalence} will always mean \textit{weak equivalence} and will refer to a model category structure.

The homotopy category of a model category $M$ is $W^{-1}M$ (see \cite[Def. $1.2.1$]{ho}),
where $W$ is the subcategory of equivalences in $M$, and it
will be denoted as $\mathrm{Ho}(M)$. We will say that two objects in a model category
$M$ are equivalent if they are isomorphic in $\mathrm{Ho}(M)$.
We say that two model categories are \textit{Quillen equivalent} if they can be connected by a finite string of
Quillen adjunctions each one being a
Quillen equivalence.

For the notions of $\mathbb{U}$-cofibrantly generated, $\mathbb{U}$-combinatorial
and $\mathbb{U}$-cellular model category, we refer to \cite{ho, hi, du2} or
to Appendix B of \cite{partI}, where the basic definitions and crucial properties are
recalled in a way that is suitable for our needs.

For a category $C$ in a universe $\mathbb{U}$, we will denote by $\mathsf{Pr}(C)$ the category of
presheaves of $\mathbb{U}$-sets
on $C$, $\mathsf{Pr}(C):=C^{\mathsf{Set}_{\mathbb{U}}^{op}}$.

\end{section}

\begin{section}{Segal categories}

As carefully explained in Leinster's survey
\cite{le}, there already exist several
definitions of higher categories, and most (if not all) of them might
be used to develop a theory of \textit{higher sites} and \textit{higher
topoi}. Because of its degree of advancement and its great
flexibility, we have chosen to work with the notion of
($1$-)\textit{Segal categories}, which might be assimilated to  $\infty$-categories
where all $i$-morphisms are invertible as soon as $i>1$.

Segal categories were first introduced by W. Dwyer, D. Kan and J. Smith
in \cite{dk3}, under the name of \textit{special bi-simplicial sets}.
With the name of \textit{$\Delta$-categories}, they were studied in more
details by R. Schw\"anzl and R. Vogt in \cite{sv}, and used in order to deal with homotopy coherence of  diagrams.
More recently, the homotopy theory of Segal categories was studied by C.
Simpson and A. Hirschowitz in \cite{sh}, and then reconsidered
in great detail by R. Pellissier in \cite{p}. They proved in particular the existence of a \textit{closed model
structure}. Using the existence of this model structure they have showed that
many of the usual categorical constructions (categories of functors, limits, colimits,
adjunctions, stacks \dots) have reasonable extensions to Segal categories.
In this Section we review briefly the main defintions in the theory of Segal categories.  \\

Recall from \cite{p,sh} that a Segal pre-category (in $\mathbb{U}$) is a functor
$$A : \Delta^{o} \longrightarrow \mathsf{SSet}_{\mathbb{U}}$$
such that $A_{0}$ is a discrete (or constant) simplicial set. The category of
Segal pre-categories in $\mathbb{U}$ will be denoted by $\mathsf{PrSeCat}_{\mathbb{U}}$ (or $\mathsf{PrSeCat}$ is the
ambient universe $\mathbb{U}$ is clear).

For any integer $n\geq 0$, one can consider the $(n+1)$-morphisms
$[0] \longrightarrow [n]$ in $\Delta$, and the associated
morphism of simplicial sets
$$A_{n} \longrightarrow A_{0}^{n+1}.$$
As the simplicial set $A_{0}$ is discrete, the existence of the above morphism implies the existence
of a natural isomorphism in $\mathsf{SSet}$
$$A_{n}\simeq \coprod_{(a_{0},\dots,a_{n})\in A_{0}^{n+1}}A_{(a_{0},\dots,a_{n})},$$
where $A_{(a_{0},\dots,a_{n})}$ is the fiber of $A_{n} \longrightarrow A_{0}^{n+1}$
at the point $(a_{0},\dots,a_{n})$.

For a Segal pre-category $A$, one should think of $A_{0}$ as the \textit{set of objects of $A$}, and
the simplicial set $A_{(a_{0},\dots,a_{n})}$ should be understood as the
\textit{space of composable morphisms} $\xymatrix{a_{0}\ar[r] & a_{1} \ar[r] & \dots \ar[r] & a_{n}}$.
In particular, $A_{(a,b)}$ is the space of morphisms from $a$ to $b$. Note for example that the
first degeneracy morphism $A_{0} \longrightarrow A_{1}$ gives natural $0$-simplices in $A_{(a,a)}$
that play the role of identities. \\

Let us fix a set $O$ in $\mathbb{U}$, and consider the subcategory $\mathsf{PrSeCat}(O)$ of Segal pre-categories
$A$ with $A_{0}=O$ and morphisms inducing the identity on $O$. Objects in $\mathsf{PrSeCat}(O)$ will be
called \textit{Segal pre-categories over the set} $O$. Note that,
obviously, any Segal pre-category $A$ can be considered as
an object in $\mathsf{PrSeCat}(A_{0})$.
A morphism $f : A \longrightarrow B$ in $\mathsf{PrSeCat}(O)$ will be called
an \textit{iso-equivalence} (resp. an \textit{iso-fibration}) if for any $(a_{0},\dots,a_{n}) \in O^{n+1}$,
the induced morphism
$$f_{(a_{0},\dots,a_{n})} : A_{(a_{0},\dots,a_{n})} \longrightarrow B_{(a_{0},\dots,a_{n})}$$
is an equivalence (resp., a fibration) of simplicial sets. There exists a simplicial closed model structure
on the category $\mathsf{PrSeCat}(O)$ where the equivalences (resp., the fibrations) are the iso-equivalences (resp., the
iso-fibrations). Indeed, the category $\mathsf{PrSeCat}(O)$ can be identified with a certain category of pointed simplicial
presheaves on a certain $\mathbb{U}$-small category,
in such a way that iso-equivalences and iso-fibrations correspond to objectwise equivalences and
objectwise fibrations. The category $\mathsf{PrSeCat}(O)$ with this model structure will be called the
\textit{objectwise model category of Segal pre-categories over $O$}. For this model structure, $\mathsf{PrSeCat}(O)$ is
furthermore a $\mathbb{U}$-combinatorial and $\mathbb{U}$-cellular model category in
the sense of \cite[Appendix]{partI} or \cite{ho, hi, du2}.

For any integer $n\geq 2$ and $0\leq i<n$, there exists a morphism
$h_{i} : [1] \longrightarrow [n]$ in $\Delta$, sending
$0$ to $i$ and $1$ to $i+1$. For any $n\geq 2$ and any
$A \in \mathsf{PrSeCat}$, the morphisms $h_{i}$ induces a morphism of simplicial
sets
$$A_{n} \longrightarrow A_{1}^{n},$$
called \textit{Segal morphism}.
This morphism breaks into several morphisms of simplicial sets
$$A_{(a_{0},\dots,a_{n})} \longrightarrow
A_{(a_{0},a_{1})}\times A_{(a_{1},a_{2})} \times\dots \times A_{(a_{n-1},a_{n})},$$
for any $(a_{0},\dots,a_{n}) \in A_{0}^{n+1}$.

For any integer $n\geq 1$ and any $(a_{0},\dots,a_{n}) \in A_{0}^{n+1}$, the functor
$$\begin{array}{cccc}
\mathsf{PrSeCat}(O) & \longrightarrow & \mathsf{SSet} \\
A & \mapsto & A_{(a_{0},\dots,a_{n})}
\end{array}$$
is co-representable by an object $h_{(a_{0},\dots,a_{n})} \in \mathsf{PrSeCat}(O)$.
The natural morphisms mentioned above,
give rise to morphisms in $\mathsf{PrSeCat}(O)$
$$h(a_{0},\dots,a_{n}) : h_{(a_{0},a_{1})}\coprod h_{(a_{1},a_{2})}\coprod \dots
\coprod h_{(a_{n-1},a_{n})} \longrightarrow
h_{(a_{0},\dots,a_{n})}.$$

The set of all morphisms $h(a_{0},\dots,a_{n})$, for all $n\geq 2$, and all $(a_{0},\dots,a_{n})\in O^{n+1}$,
belongs to $\mathbb{U}$. As the objectwise model category $\mathsf{PrSeCat}(O)$ is
a $\mathbb{U}$-cellular and a $\mathbb{U}$-combinatorial model category, the following definition makes sense.

\begin{df}\label{da1}
Let $O$ be a $\mathbb{U}$-set. The \emph{model category of Segal pre-categories
over} $O$ in $\mathbb{U}$ is the left Bousfield localization
of the objectwise model category $\mathsf{PrSeCat}(O)$ with respect to the
morphisms $h(a_{0},\dots,a_{n})$, for all $n\geq 2$, and all $(a_{0},\dots,a_{n})\in O^{n+1}$.
This model category will simply be denoted by $\mathsf{PrSeCat}(O)$.
\end{df}

\textbf{Notations.} The equivalences in the model
category structure $\mathsf{PrSeCat}(O)$ of Definition \ref{da1} will simply be
called \textit{equivalences}. The homotopy category $\mathrm{Ho}(\mathsf{PrSeCat}(O))$
will always refer to the model structure $\mathsf{PrSeCat}(O)$ of the above definition. while the homotopy category of
the objectwise model structure will be denoted by $\mathrm{Ho}^{iso}(\mathsf{PrSeCat}(O))$. \\

It is clear that an object $A \in \mathsf{PrSeCat}(O)$ is fibrant if and only if it satisfies the following two
conditions

\begin{enumerate}
\item For any $n\geq 1$ and any $(a_{0},\dots,a_{n})\in O^{n+1}$, the simplicial set $A_{(a_{0},\dots,a_{n})}$
is fibrant.

\item For any $n\geq 2$ and any $(a_{0},\dots,a_{n})\in O^{n+1}$, the natural morphism
$$A_{(a_{0},\dots,a_{n})} \longrightarrow A_{(a_{0},a_{1})}\times A_{(a_{1},a_{2})} \times \dots A_{(a_{n-1},a_{n})}$$
is an equivalence of simplicial sets.
\end{enumerate}

The general theory of left Bousfield localization (see \cite{hi}) tells us that
$\mathrm{Ho}(\mathsf{PrSeCat}(O))$ is naturally equivalent
to the full subcategory of $\mathrm{Ho}^{iso}(\mathsf{PrSeCat}(O))$
consisting of objects satisfying condition $(2)$ above. Furthermore, the natural inclusion
$\mathrm{Ho}(\mathsf{PrSeCat}(O)) \longrightarrow \mathrm{Ho}^{iso}(\mathsf{PrSeCat}(O))$ possesses a left adjoint
$$SeCat : \mathrm{Ho}^{iso}(\mathsf{PrSeCat}(O)) \longrightarrow \mathrm{Ho}(\mathsf{PrSeCat}(O))$$
which is the left derived functor of the identity functor on $\mathsf{PrSeCat}(O)$.

\begin{df}\label{da2}
\begin{enumerate}
\item A Segal pre-category $A$ is a \emph{Segal category} if for any
$n\geq 2$ and any $(a_{0},\dots,a_{n})\in A_{0}^{n+1}$, the natural morphism
$$A_{(a_{0},\dots,a_{n})} \longrightarrow
A_{(a_{0},a_{1})}\times A_{(a_{1},a_{2})} \times \dots A_{(a_{n-1},a_{n})}$$
is an equivalence of simplicial sets.

\item The Segal category \emph{associated} to a Segal pre-category
$A$ over $A_{0}$, is $SeCat(A)\in \mathrm{Ho}(\mathsf{PrSeCat}(A_{0}))\subset \mathrm{Ho}^{iso}(\mathsf{PrSeCat}(A_{0}))$.

\end{enumerate}
\end{df}

Any $\mathbb{U}$-small category $C$ might be considered as a Segal category, with the same set of objects $C_{0}$ and
such that
$$C_{(a_{0},\dots,a_{n})}:=Hom_{C}(a_{0},a_{1})\times Hom_{C}(a_{1},a_{2}) \times \dots Hom_{C}(a_{n-1},a_{n})$$
for any $(a_{0},\dots,a_{n})\in C_{0}^{n+1}$. In other words,
the simplicial object $C : \Delta^{o} \longrightarrow \mathsf{SSet}$
is simply the \textit{nerve} of the category $C$, which obviously satisfied condition $(1)$ of definition \ref{da2}
above.
This allows us to view the category of $\mathbb{U}$-small categories
as a full subcategory of $\mathsf{PrSeCat}$.

If $T$ is an $S$-category (in $\mathbb{U}$), i.e. a categry enriched over simplicial sets, we can consider it as a Segal
pre-category $T : \Delta^{o} \longrightarrow \mathsf{SSet}$
by defining
$$T_{0}:=Ob(T) \qquad T_{n}:=\coprod_{(a_{0},\dots,a_{n})\in T_{0}^{n+1}}
\underline{Hom}_{T}(a_{0},a_{1})\times \dots \times
\underline{Hom}_{T}(a_{n-1},a_{n}),$$
the simplicial and degeneracies morphism being induced by the compositions
and identities in $T$. As a Segal pre-category, $T$
has the special property that the natural morphisms
$$T_{(a_{0},\dots,a_{n})} \longrightarrow T_{(a_{0},a_{1})}\times \dots \times T_{(a_{n-1},a_{n})}$$
are all \textit{isomorphisms} of simplicial sets. In particular, $T$ is
always a Segal category. This allows us to view the category $S-\mathsf{Cat}$ of $S$-categories as
the full subcategory of $\mathsf{PrSeCat}$ consisting of objects $A$ such that all morphisms
$$A_{(a_{0},\dots,a_{n})} \longrightarrow A_{(a_{0},a_{1})}\times \dots
\times A_{(a_{n-1},a_{n})}$$
are isomorphisms of simplicial sets.

More generally, a Segal category $A$ might be seen as a
\textit{weak category in $\mathsf{SSet}$}. Indeed, for any two objects
$(a,b) \in A_{0}^{2}$, the simplicial set $A_{(a,b)}$ can be considered as the space of morphisms from $a$ to $b$.
The composition of morphisms is given by a diagram
$$\xymatrix{
A_{(a,b,c)} \ar[d]_-{\sim} \ar[r] & A_{(a,c)} \\
A_{(a,b)}\times A_{(b,c)} & }$$
where the horizontal morphism is induced by the morphism $[1] \rightarrow [2]$ in $\Delta$ which sends
$0$ to $0$ and $1$ to $2$. The vertical morphism being an equivalence of simplicial sets, this diagram gives a
\textit{weak composition morphism}, i.e. a morphism
$A_{(a,b)}\times A_{(b,c)} \longrightarrow A_{(a,c)}$, which is well defined
only up to equivalence. From this point of view, the higher simplicial identities of $A$ should
be seen as providing the associativity and higher coherency
laws for this composition.

If $A$ is a Segal category, one can define the \textit{homotopy category}
$\mathrm{Ho}(A)$, whose set of objects is $A_{0}$
and whose set of morphisms from $a$ to $b$ is $\pi_{0}(A_{(a,b)})$.
The composition of morphisms in $\mathrm{Ho}(A)$
is given by the following diagram
$$\xymatrix{
\pi_{0}(A_{(a,b,c)}) \ar[r] \ar[d]_-{\sim}  & \pi_{0}(A_{(a,c)}) \\
\pi_{0}(A_{(a,b)})\times \pi_{0}(A_{(b,c)}), & }$$
where the horizontal morphism is induced by the morphism $[1] \rightarrow [2]$ in $\Delta$ which sends
$0$ to $0$ and $1$ to $2$.

Any Segal pre-category $A$ can be considered as an object in $\mathsf{PrSeCat}(A_{0})$,
where $A_{0}$ is the set of objets of $A$.
Applying the functor $SeCat$, one finds a natural morphism in $\mathrm{Ho}^{iso}(\mathsf{PrSeCat}(A_{0}))$
$$A \longrightarrow SeCat(A).$$
This morphism has to be thought of as a (Segal) \textit{categorification} of $A$.
This way, one can define the homotopy category of the Segal pre-category $A$ to be
$\mathrm{Ho}(SeCat(A))$. In the case $A$ is already a Segal category, the natural
morphism $A \longrightarrow SeCat(A)$ is an iso-equivalence, and therefore induces
an isomorphism $\mathrm{Ho}(A) \simeq \mathrm{Ho}(SeCat(A))$. This shows that the homotopy category
functor $A \mapsto \mathrm{Ho}(A)$ is well defined.

An object in a Segal category $A$ is an element in $A_{0}$. Most of the time
they will be considered up to equivalences, or in other words as
isomorphism classes of objects in $\mathrm{Ho}(A)$. In the same way, a morphism
$a \rightarrow b$ in $A$ is a $0$-simplex in $A_{(a,b)}$, mots often
considered as a morphism in $\mathrm{Ho}(A)$ (i.e. up to homotopy).

A morphism in a Segal category $A$ (i.e. an element of $A_{a,b}$ for some $(a,b)\in A^{2}_{0}$)
is called an \textit{equivalence} if its image in $\mathrm{Ho}(A)$ is an isomorphism.
There exists a sub-Segal category $A^{int}$ of $A$, consisting of all objects and
equivalences between them. The Segal category $A^{int}$ is a Segal groupoid (i.e.
all its morphisms are equivalences) and is actually the maximal sub-Segal groupoid
of $A$ (see \cite[\S 2]{sh}).

Let $f : A \longrightarrow B$ be a morphism in $\mathsf{PrSeCat}$.
This morphism induces a morphism of sets $f : A_{0} \longrightarrow B_{0}$, as
well as a morphism in $\mathsf{PrSeCat}(A_{0})$
$$A \longrightarrow f^{*}(B),$$
where $f^{*}(B)$ is defined by the formula
$$f^{*}(B)_{(a_{0},\dots,a_{n})}:=B_{(f(a_{0}),\dots,f(a_{n})} \qquad \; (a_{0},\dots,a_{n})\in A_{0}^{n+1}.$$

Composing with the morphism $B \longrightarrow SeCat(B)$, one gets a morphism
$$A \longrightarrow f^{*}(SeCat(B))$$
which is well defined in $\mathrm{Ho}^{iso}(\mathsf{PrSeCat}(A_{0}))$. The Segal pre-category $f^{*}(SeCat(B))$ is
clearly a Segal category, and therefore there exists a unique factorization in $\mathrm{Ho}^{iso}(\mathsf{PrSeCat}(A_{0}))$
$$\xymatrix{
A \ar[d] \ar[r] & f^{*}(SeCat(B)) \\
SeCat(A). \ar[ru] & }$$
In particular, there are natural morphisms in $\mathrm{Ho}(\mathsf{SSet})$
$$SeCat(A)_{(a,b)} \longrightarrow SeCat(B)_{(f(a),f(b))}$$
for any $(a,b) \in A_{0}^{2}$. The induced morphism on the connected components, together with
the map $A_{0} \longrightarrow B_{0}$, give rise to a
well defined functor
$$\mathrm{Ho}(f) : \mathrm{Ho}(SeCat(A)) \longrightarrow \mathrm{Ho}(SeCat(B)).$$

\begin{df}\label{da3}
Let $f : A \longrightarrow B$ be a morphism of Segal pre-categories.
\begin{itemize}
\item The morphism $f$ is \emph{fully faithful}
if for any $(a,b) \in A_{0}^{2}$, the induced morphism
$$SeCat(A)_{(a,b)} \longrightarrow SeCat(B)_{(f(a),f(b))}$$
is an isomorphism in $\mathrm{Ho}(\mathsf{SSet})$.

\item The morphism $f$ is \emph{essentially surjective} if the induced functor
$$\mathrm{Ho}(SeCat(A)) \longrightarrow \mathrm{Ho}(SeCat(B))$$
is essentially surjective.

\item The morphism $f$ is an \emph{equivalence} if it is fully faithful and essentially surjective.
\end{itemize}
\end{df}

It is important to notice that when $f : A \longrightarrow B$ is a morphism between Segal categories, then
$A$ and $B$ are iso-equivalent to $SeCat(A)$ and $SeCat(B)$, and therefore the previous definition simplifies. \\

\begin{thm}\emph{(\cite[Thm. $6.4.4$]{p})}
There exists a closed model structure on the category
$\mathsf{PrSeCat}$, of Segal pre-categories (in the universe $\mathbb{U}$),
where the equivalences are those of Definition \ref{da3} and the cofibrations are the monomorphisms. This model category
is furthermore simplicial, $\mathbb{U}$-cellular and $\mathbb{U}$-combinatorial.
\end{thm}

For any $A$ and $B$ in $\mathsf{PrSeCat}$, the natural morphism in $\mathrm{Ho}^{iso}(\mathsf{PrSeCat})$
$$SeCat(A\times B) \longrightarrow SeCat(A)\times SeCat(B)$$
is an isomorphism. This is the \textit{product formula} of \cite[Thm. $5.5.20$]{p}. From this, one deduces formally the
following very important corollary.

\begin{cor}\emph{(\cite[Thm. $5.5.20$, Thm. $6.4.4$]{p})}\label{ca1}
The model category $\mathsf{PrSeCat}$ is internal (i.e. is a symmetric monoidal model category
with the direct product as the monoidal product).
\end{cor}

As explained in \cite[Thm. $4.3.2$]{ho}, the previous corollary implies that the homotopy category
$\mathrm{Ho}(\mathsf{PrSeCat})$ is cartesian closed; the correponding internal $Hom$-objects will be denoted by
$\mathbb{R}\underline{Hom}(A,B) \in \mathrm{Ho}(\mathsf{PrSeCat})$,
for $A$ and $B$ in $\mathrm{Ho}(\mathsf{PrSeCat})$. Recall that
$\mathbb{R}\underline{Hom}(A,B)$ is naturally isomorphic to $\underline{Hom}(A,RB)$, where
$\underline{Hom}$ denotes the internal $Hom$'s in the category $\mathsf{PrSeCat}$ and $RB$ is a fibrant model of $B$
in $\mathsf{PrSeCat}$. For any $A$, $B$ and $C$ in $\mathrm{Ho}(\mathsf{PrSeCat})$,
one has the derived  adjunction formula
$$\mathbb{R}\underline{Hom}(A,\mathbb{R}\underline{Hom}(B,C))\simeq \mathbb{R}\underline{Hom}(A\times B,C).$$

The existence of these derived internal $Hom$'s is of fundamental importance as it allows one to
develop the theory of Segal categories in a very similar fashion as usual category theory.
A first example is the following definition. For any category $C$ together with a subcategory $S$,
we denote by $L(C,S)$ the Dwyer-Kan simplicial localization of $C$ with respect to $S$
(see \cite{dk1}); $L(C,S)$ is an $S$-category hence a Segal category. When $C$ is a model
category we will always take $S$ to be the subcategory of weak equivalences and we will simply write $LC$ for $L(C,S)$.

\begin{df}\label{da4}
\begin{enumerate}
\item The Segal category of $\mathbb{U}$-small simplicial sets is
defined to be
$$\mathsf{Top}:=L(\mathsf{SSet}_{\mathbb{U}}).$$
\item
Let $T$ be a $\mathbb{U}$-small Segal category. The \emph{Segal category of pre-stacks over}
$T$ is defined to be
$$\hat{T}:=\mathbb{R}\underline{Hom}(T^{op},\mathsf{Top}).$$
\end{enumerate}
\end{df}

Note that, as usual, if $T$ is a $\mathbb{U}$-small Segal category, then $\hat{T}$ is
a $\mathbb{V}$-small Segal category for $\mathbb{U}\in \mathbb{V}$.

\end{section}

\begin{section}{Segal topoi}

The existence of internal $Hom$'s mentioned above, and more
generally of the model strcuture on $\mathsf{PrSeCat}$, allows one to
extend to Segal categories most of the basic constructions in category theory (for some of these constructions,
see \cite{sh,s3}). In this
paragraph we will first recall some of the basic extensions and then use them to describe a notion of
\textit{Segal topos}, analogous to the usual notion of topos
in the context of Segal categories. We will state the results without proofs.

\begin{subsection}{Adjunctions}

Let us start by the definition of \textit{adjunction} between Segal categories, as introduced in \cite[\S $8$]{sh}.
Let $f : A \longrightarrow B$ be a morphism between fibrant Segal categories.
We will say that $f$ \textit{has a left adjoint}, if there
exists a morphism $g : B \longrightarrow A$ and an element $u \in \underline{Hom}(A,A)_{gf,Id}$ (i.e. a natural
transformation $gf \Rightarrow Id$) such that for any $a \in A_{0}$ and $b\in B_{0}$, the natural morphism
(well defined in $\mathrm{Ho}(\mathsf{SSet})$)
$$\xymatrix{
B_{(f(a),b)} \ar[r]^-{g_{*}} & A_{(gf(a),g(b))} \ar[r]^-{u_{*}} & A_{(a,g(b))}}$$
is an isomorphism in $\mathrm{Ho}(\mathsf{SSet})$.

One can check that if $f$ has a left adjoint, then the pair $(g,u)$ is unique up to equivalence. More
generally one can show that the Segal category of left adjoints to $f$ (i.e. of all pairs $(g,u)$ as above)
is contractible. This justifies the terminology
\textit{$g$ is the left adjoint of $f$, and $u$ is its unit}, for any such pair $(g,u)$.

\begin{df}
A morphism $f : A \longrightarrow B$ in $\mathrm{Ho}(\mathsf{PrSeCat})$ \emph{has a left adjoint}
if for any $A'$ and $B'$ fibrant models of $A$ and $B$,
one of its representative $f : A' \longrightarrow B'$ in $\mathsf{PrSeCat}$ has a left adjoint in the sense
defined above.
\end{df}

In a dual way one defines the notion of right adjoints.

\end{subsection}

\begin{subsection}{Limits}

Let $I$ and $A$ be two Segal categories. The natural projection $I\times A \longrightarrow A$ yields by adjunction
a well defined morphism in $\mathrm{Ho}(\mathsf{PrSeCat})$
$$c : A \longrightarrow \mathbb{R}\underline{Hom}(I,A).$$
The morphism $c$ is the \textit{constant diagram functor}, and sends an object $a \in A_{0}$ to the constant
morphism $I \rightarrow \{a\} \subset A$.

\begin{df}
\begin{enumerate}
\item A Segal category $A$ has \emph{limits} (resp. \emph{colimits}) \emph{along a category} $I$ if
the natural morphism $c : A \longrightarrow \mathbb{R}\underline{Hom}(I,A)$ has a right (resp. left) adjoint.
The right adjoint (resp. left adjoint) is then denoted by
$$\mathrm{Lim}_{I} : \mathbb{R}\underline{Hom}(I,A) \longrightarrow A \qquad
\left(resp. \; \mathrm{Colim}_{I} \mathbb{R}\underline{Hom}(I,A) \longrightarrow A \right).$$

\item
Let $\mathbb{U}\in \mathbb{V}$ be two universes, and
$A$ a Segal category in $\mathbb{V}$, $A \in \mathrm{Ho}(\mathsf{PrSeCat}_{\mathbb{V}})$.
The Segal category $A$ has $\mathbb{U}$\emph{-limits} (resp. $\mathbb{U}$\emph{-colimits}),
if it has limits (resp. colimits)
along any $\mathbb{U}$-small category $I$.

\item A Segal category $A$ has \emph{finite limits} if it has
limits along the categories $(\xymatrix{2 \ar[r] & 0 & \ar[l] 1})$ and $\; \emptyset$. A Segal category has
\emph{finite colimits} if its opposite Segal category $A^{op}$ has finite limits.

\end{enumerate}
\end{df}

When concerned with limits and colimits in Segal categories we will use the same
notations as usual. For example, coproducts will be denoted
by $\coprod$, fibered products by $x\times_{z}y$, \dots. \\

\begin{rmk} \emph{One can show that any Segal category $A$ with finite limits in the sense above has also
limits along all categories $I$ whose nerve is a finite simplicial set.
This justifies the terminology of having \textit{finite limits}.}
\end{rmk}

The following proposition is a formal consequence of the definitions and
the adjunction formula for internal Hom's.

\begin{prop}
For a $\mathbb{U}$-small Segal category $T$, the Segal category $\hat{T}$ has $\mathbb{U}$-limits and
$\mathbb{U}$-colimits.
\end{prop}

Let $f : A\longrightarrow B$ be a morphism of Segal categories, and $I$ any category such that
$A$ and $B$ have limits along $I$. The universal property of adjunction implies that for any
$x_{*}\in \mathbb{R}\underline{Hom}(I,A)$, there exists a well defined and natural morphism
in $\mathrm{Ho}(B)$
$$f(\mathrm{Lim}_{I}x_{i}) \longrightarrow \mathrm{Lim}_{I}f(x_{i}).$$
We will say that $f$ \textit{preserves limits along} $I$ if for any $x_{*} \in \mathbb{R}\underline{Hom}(A,B)$,
the induced morphism $f(\mathrm{Lim}_{I}x_{i})
\longrightarrow \mathrm{Lim}_{I}f(x_{i})$ is an isomorphism in $\mathrm{Ho}(B)$.
One shows, for example, that a morphism having a left adjoint always preserves limits.

\begin{df}\label{dlexSeg}
Let $A$ and $B$ be two Segal categories with finite limits.
\begin{enumerate}
\item
A morphism $f : A \longrightarrow B$ in $\mathrm{Ho}(\mathsf{PrSeCat})$ is \emph{left exact}
if it preserves limits along $\xymatrix{2 \ar[r] & 0 & \ar[l] 1}$ and
$\emptyset$.

\item
The Segal category $A$ is a \emph{left exact localization of} $B$ if there exists
a morphism $i : A \longrightarrow B$ in $\mathrm{Ho}(\mathsf{PrSeCat})$, which is fully
faithful and possesses a left adjoint which
is left exact.
\end{enumerate}
\end{df}

In a dual way one defines the notion of \textit{colimits preserving} morphism,
\textit{right exact} morphism and \textit{right exact localization}.

\end{subsection}

\begin{subsection}{Topologies and Segal topoi}

\begin{df}\label{dtop}
A \emph{Segal topology on a Segal category} $T$ is a Grothendieck
topology on the homotopy category $\mathrm{Ho}(T)$.
A Segal category together with a Segal topology is called a \emph{Segal site}.
\end{df}

The functor $\pi_{0}$, sending a simplicial set to its set of connected
components, induces a morphism of Segal categories
$\pi_{0} : \mathsf{Top} \longrightarrow \mathsf{Set}$. By composition, one gets for any
$\mathbb{U}$-small Segal category $T$, a morphism
$$\pi^{pr}_{0} : \hat{T} \longrightarrow \mathbb{R}\underline{Hom}
(T,\mathsf{Set})\simeq \mathbb{R}\underline{Hom}(\mathrm{Ho}(T),\mathsf{Set})\simeq
\mathsf{Pr}(\mathrm{Ho}(T)),$$
where $\mathsf{Pr}(\mathrm{Ho}(T))$ is the category of presheaves of sets
on the (usual) category $\mathrm{Ho}(T)$.

In the same way, for any simplicial set $K$, the exponential functor $(-)^{K} : \mathsf{SSet}
\longrightarrow \mathsf{SSet}$, restricted to
the full subcategory of fibrant objects, induces a morphism of Segal categories
$$(-)^{\mathbb{R}K} : \mathsf{Top} \longrightarrow \mathsf{Top}.$$
By composition, this induces a morphism of Segal categories of pre-stacks
$$(-)^{\mathbb{R}K} : \hat{T} \longrightarrow \hat{T}.$$

\begin{df}\label{dstacks}
Let $(T,\tau)$ be a $\mathbb{U}$-small Segal site.
\begin{enumerate}
\item A morphism $f : F \longrightarrow G$ in $\hat{T}$ is a $\tau$\emph{-local equivalence} if
for any integer $n\geq 0$, the induced morphism of presheaves of sets on $\mathrm{Ho}(T)$
$$\pi^{pr}_{0}(F^{\mathbb{R}\partial \Delta^{n}}) \longrightarrow
\pi^{pr}_{0}(F^{\mathbb{R}\partial \Delta^{n}}\times_{G^{\mathbb{R}\partial \Delta^{n}}}G^{\mathbb{R}\Delta^{n}})$$
induces an epimorphism of sheaves on the site $(\mathrm{Ho}(T),\tau)$.

\item An object $F \in (\hat{T})_{0}$ is called a \emph{stack} for the topology $\tau$, if for any
$\tau$-local equivalence $f : G \longrightarrow H$ in $\hat{T}$, the induced morphism
$$f^{*} : Hom_{\mathrm{Ho}(\hat{T})}(H,F) \longrightarrow Hom_{\mathrm{Ho}(\hat{T})}(G,F)$$
is bijective.

\item The \emph{Segal category of stacks} on the Segal site $(T,\tau)$ is
the full sub-Segal category $T^{\sim,\tau}$ of $\hat{T}$, consisting of stacks.

\end{enumerate}
\end{df}

As mentioned in the Introduction of \cite{partI}, one possible definition of a usual Grothendieck
topos is as a full subcategory of a category of presheaves such that
the inclusion functor possesses a left adjoint which is left exact. The following definition
is an analog of this for Segal categories.

\begin{df}\label{dsegtop}
A $\mathbb{U}$\emph{-Segal topos} is a Segal category which is a left exact localization of
$\hat{T}$, for a $\mathbb{U}$-small Segal category $T$.
\end{df}

\begin{rmk}
\emph{
\begin{itemize}
\item
When $T$ is a $\mathbb{U}$-small Segal category, $\hat{T}$ is only $\mathbb{V}$-small for
$\mathbb{U}\in \mathbb{V}$. Therefore, a $\mathbb{U}$-Segal topos does not belong to $\mathbb{U}$ but only
to $\mathbb{V}$.
\item
By definition of left exact localizations in \ref{dlexSeg}, a $\mathbb{U}$-Segal
category $A$ is a Segal topos if and only if
there exists a $\mathbb{U}$-small Segal category $T$ and a morphism
$i : A \longrightarrow \hat{T}$, which is
fully faithful and has a left exact left adjoint. The morphism $i$
and the Segal category $T$ are however \textit{not} part of the data.
\end{itemize}
}
\end{rmk}

As a direct consequences of the definition one has the following result.

\begin{prop}
A $\mathbb{U}$-Segal topos has $\mathbb{U}$-limits and $\mathbb{U}$-colimits.
\end{prop}

In a similar way as in \cite[Def. 3.8.2]{partI}, we define the notion of truncated objects in a Segal category $A$. An
object $a \in A$ is $n$\textit{-truncated}, if for any $x \in A$
the simplicial set $A_{(x,a)}$ is $n$-truncated. An object
will be simply called \textit{truncated} if it is $n$-truncated for some $n$.

\begin{df}
A Segal category $A$ is
$t$\emph{-complete} if truncated objects detect isomorphisms in
$\mathrm{Ho}(A)$ (i.e. a morphism $u : x \rightarrow y$ is an isomorphism
in $\mathrm{Ho}(A)$ if and only if $u^{*} :
[y,a]_{\mathrm{Ho}(A)}\rightarrow  [x,a]_{\mathrm{Ho}(A)}$
is bijective for any truncated object $a \in \mathrm{Ho}(A)$).
\end{df}

The fundamental example of a $t$-complete Segal category is the following.

\begin{prop}
\begin{enumerate}
\item
Let $(T,\tau)$ be $\mathbb{U}$-small Segal site. The natural inclusion functor
$T^{\sim,\tau} \longrightarrow \hat{T}$ has a left exact left adjoint.

\item
The Segal category $T^{\sim,\tau}$ of stacks on $(T,\tau)$
is a $t$-complete $\mathbb{U}$-Segal topos.
\end{enumerate}
\end{prop}

\textit{Proof:} It is a consequence of the corresponding fact
in the context of model categories (see \cite[Proposition 3.4.10 (2)]{partI}
and \cite[Theorem 3.8.3]{partI}) and the
comparison theorem \ref{ccomp}. \hfill $\Box$ \\

We are now ready to state the analog of Theorem 3.8.3 of \cite{partI} for Segal categories.
We fix a $\mathbb{U}$-small Segal category
$T$. Any topology $\tau$ on $T$ gives a full subcategory of stacks $T^{\sim,\tau}$,
which by the previous proposition, is
a $t$-complete left exact localization of $\hat{T}$.  The proof of the following theorem can
be established, with some work, using the same statement for $S$-categories (proved in \cite[Theorem 3.8.3]{partI}), and
the comparison result of Corollary \ref{ccomp}.

\begin{thm}\label{cmp}
With the above notations, the rule $\tau \longmapsto T^{\sim,\tau}$ establishes
a bijection between Segal topologies on $T$ and
(equivalence classes) of $t$-complete left exact localizations of $\hat{T}$.
\end{thm}

\begin{rmk}\label{hyper}\emph{The hypothesis of $t$-completeness in Theorem \ref{cmp} might appear unnatural,
and it would be interesting to understand whether there exists a kind of ``topologies''
on $T$ which are in bijection with \textit{arbitrary} left exact localizations of $\hat{T}$.
As explained in \cite{partI} Remark 3.8.7 (3), one way to proceed would be to introduce a  \textit{hyper-topology} on a
Segal category, a notion suggested to us by some remarks of
V. Hinich, A. Joyal and C. Simpson.
A hyper-topology on a (Segal) category would consist in specifying directly the hypercovers
(and not only the coverings, like in the case of a topology). The data of these hypercovers
should satisfy appropriate conditions ensuring that the ``correponding'' left Bousfield
localization is indeed exact. Then, it seems reasonable that Theorem \ref{cmp} can be
generalized to a bijective correspondence between hyper-topologies
on $T$ and arbitrary left exact Bousfield localizations of $\hat{T}$.}

\emph{Theorem \ref{cmp} suggests also a way to think of \textit{higher topologies} on $n$-Segal
categories (and of \textit{higher topoi}) for $n \geq 1$ as appropriate \textit{left exact
localizations}. We address the reader to \cite[Remark 3.8.7 (4)]{partI} for a brief discussion on this point.}
\end{rmk}
\end{subsection}

\begin{subsection}{Geometric morphisms}

The main advantage of Segal categories with respect to $S$-categories is the existence of a reasonable theory of
\textit{internal} $Hom$'s (see Corollary \ref{ca1}).
For the purpose of topos theory, this will allow us to define the notion of
\textit{geometric morphisms} between Segal topoi, and more generally of the
\textit{Segal category of geometric morphisms} between
two Segal topoi.

\begin{df}\label{dgeom}
Let $A$ and $B$ be two $\mathbb{U}$-Segal topoi.
\begin{enumerate}
\item
A morphism $f : B \longrightarrow A$ in $\mathrm{Ho}(\mathsf{PrSeCat})$ is
called \emph{geometric} if it is left exact and has a right adjoint.

\item
The full sub-Segal category of $\mathbb{R}\underline{Hom}(B,A)$ consisting of geometric morphisms will be
denoted by $\mathbb{R}\underline{Hom}^{geom}(A,B)$.

\end{enumerate}
\end{df}

\begin{rmk}\begin{enumerate}
\item \emph{As one uses the internal $Hom$'s for Segal categories in order to
define the $2$-Segal category of Segal categories (see
\cite[\S $2$]{sh}),
one can use the notion of geometric morphisms between Segal topoi to
define the $2$\textit{-Segal category of Segal topoi}.
More precisely, one defines
the $2$-Segal category of $\mathbb{U}$-Segal topoi $\mathsf{SeT}_{\mathbb{U}}$ as
follows. Its objects are the fibrant Segal categories in $\mathbb{V}$ which are $\mathbb{U}$-Segal topoi.
For $A$ and $B$ two
$\mathbb{U}$-Segal topoi, the $1$-Segal category of morphisms
is $\underline{Hom}^{geom}(A,B)$, the full sub-Segal category of
$\underline{Hom}(B,A)$ consisting of geometric morphisms.
Note that $\mathsf{SeT}_{\mathbb{U}}$ is an element of $\mathbb{W}$, for a
universe $\mathbb{W}$ such that $\mathbb{U}\in \mathbb{V} \in \mathbb{W}$.}

\item \emph{In part $(1)$ of Definition \ref{dgeom}, the right adjoint is not part of the data of a geometric
morphism. Strictly speaking, this might differ from the original definition of geometric morphisms, but the
uniqueness of adjoints implies that the two notions are in fact equivalent (i.e. give rise to equivalent
$2$-categories of topoi).}

\item \emph{There is also a notion of \textit{essential} geometric
morphism of topoi as described for example, in \cite[p. 360]{mm}. In the same
way one can define the refined notion of \textit{essential geometric morphisms of Segal topoi}.
They form a full sub-Segal category
of the Segal category of geometric morphisms.}

\end{enumerate}

\end{rmk}

The standard example of a geometric morphism of Segal topoi is given by
continuous morphisms of Segal sites. Let $(T,\tau)$ and
$(T',\tau')$ be two $\mathbb{U}$-small Segal sites, and $f : T \longrightarrow T'$
be a morphism in $\mathrm{Ho}(\mathsf{PrSeCat})$. The morphism $f$ induces a
morphism in $\mathrm{Ho}(\mathsf{PrSeCat})$ between the corresponding Segal categories of pre-stacks
$$f^{*} : \hat{T'} \longrightarrow \hat{T}.$$
We say that the morphism $f$ is \textit{continuous} if it preserves the full
sub-Segal category of stacks. In this case, it induces a well
defined morphim in $\mathrm{Ho}(\mathsf{PrSeCat})$
$$f^{*} : (T')^{\sim,\tau'} \longrightarrow T^{\sim,\tau}.$$
One can show that this morphism is left exact and possesses a left and a
right adjoint, $f_{!}$ and $f_{*}$. In particular,
$f^{*}$ defines a geometric morphism $T^{\sim,\tau} \longrightarrow (T')^{\sim,\tau}$.

\end{subsection}

\end{section}

\begin{section}{Comparison between model topoi and Segal topoi}

\begin{subsection}{$S$-categories vs. Segal categories}\label{svssegal}

As we already observed, any $S$-category (i.e. simplicially enriched category) is in a natural way a Segal category.
Furthermore, it is clear that a
morphism $T \longrightarrow T'$ between $S$-categories
is an equivalence if and only if it is an equivalence in the model
category $\mathsf{PrSeCat}$. Therefore, the inclusion functor
$j : S-\mathsf{Cat} \longrightarrow \mathsf{PrSeCat}$ induces a well defined functor on the level of homotopy categories
$$j : \mathrm{Ho}(S-\mathsf{Cat}) \longrightarrow \mathrm{Ho}(\mathsf{PrSeCat}).$$
This functor is known to be an equivalence of categories (see \cite[p. $7$]{s3}).
Since we know that $\mathrm{Ho}(\mathsf{PrSeCat})$ is
cartesian closed, then so is $\mathrm{Ho}(S-\mathsf{Cat})$. \\

Let $C$ be a $\mathbb{U}$-small category, $S \subset C$ a subcategory
(that we may suppose to contain all the isomorphisms) and $L(C,S)$ its simplicial localization,
considered as a Segal category and therefore as an object in $\mathrm{Ho}(\mathsf{PrSeCat})$.
It comes equipped with a natural localization morphism
$L : C \longrightarrow \mathrm{L}(C,S)$ in $\mathrm{Ho}(\mathsf{PrSeCat})$
defined by the diagram $$\xymatrix{L(C,\mathrm{iso})\ar[r] \ar[d]_{\sim} & L(C,S)\\
C & }$$ where the horizontal map is induced by the identity functor on $C$
and the vertical equivalence in $\mathsf{PrSeCat}$ is adjoint to the
(usual) equivalence of categories $\pi_{0}L(C,\mathrm{iso})\simeq C$.

The following proposition says that
$C \longrightarrow \mathrm{L}(C,S)$ is the
universal construction
in $\mathrm{Ho}(\mathsf{PrSeCat})$ which \textit{formally inverts} the morphisms in $S$.

\begin{prop}\emph{(\cite[Prop. $8.6$, Prop. $8.7$]{sh})}\label{uni}
For any $A \in \mathrm{Ho}(\mathsf{PrSeCat})$, the natural morphism
$$L^{*} : \mathbb{R}\underline{Hom}(\mathrm{L}(C,S),A) \longrightarrow
\mathbb{R}\underline{Hom}(C,A)$$
is fully faithful. Its essential image consists of those morphisms $C \rightarrow A$ mapping elements of $S$ to
equivalences in $A$.
\end{prop}

Proposition \ref{uni} has the following important corollary.

\begin{cor}
Let $(C,S)$ (resp. $(D,T)$) be a $\mathbb{U}$-small category with a subcategory
$S \subset C$ (resp. $T \subset D$). Then, the natural morphism in
$\mathrm{Ho}(\mathsf{PrSeCat})$
$$L(C\times D,S\times T) \longrightarrow L(C,S)\times L(D,T)$$
is an isomorphism.
\end{cor}

\end{subsection}

\begin{subsection}{The strictification theorem}

Let $T$ be a $\mathbb{U}$-small $S$-category, and $M$ be a simplicial model category which will be
assumed to be $\mathbb{U}$-cofibrantly generated. The category of simplicial functors from $T$ to $M$ is denoted by
$M^{T}$, and is endowed with its projective objectwise model structure:
equivalences and fibrations are defined objectwise. We consider the
\textit{evaluation functor}
$$M^{T}\times T \longrightarrow M,$$
as a morphism in $\mathsf{PrSeCat}$. This functor clearly sends objectwise
equivalences in $M^{T}\times T$ to equivalences in $M$,
and therefore induces a well defined morphism in $\mathrm{Ho}(\mathsf{PrSeCat})$
between the corresponding simplicial localizations along equivalences
$$L(M^{T}\times T)\simeq L(M^{T})\times T \longrightarrow LM,$$
which yields, by adjunction, a morphism
$$L(M^{T}) \longrightarrow \mathbb{R}\underline{Hom}(T,LM).$$

The \textit{strictification theorem} is the following statement.
It says that the Segal category of functors from $T$ to $LM$ can
be computed using the model category $M^{T}$ of $T$-diagrams in $M$.

\begin{thm}\emph{(Strictification theorem)}\label{tstrict}
Under the previous hypotheses and notations, the natural morphism
$$L(M^{T}) \longrightarrow \mathbb{R}\underline{Hom}(T,LM)$$
is an isomorphism.
\end{thm}

\textit{Sketch of proof:} The theorem is proved in \cite[Theorem 18.6]{sh} when $T$ is
a Reedy category. In the general case, one can find a Reedy category $C$ and a subcategory $S\subset C$
together with an equivalence $L(C,S) \simeq T$. The theorem now follows from the case of
a Reedy category, Proposition \ref{uni} and its analog for model categories
(see \cite[Theorem 2.3.5]{partI} for references). \hfill $\Box$ \\

This theorem is the central result needed to compare the construction of the category of stacks on an
$S$-site (\cite[\S 3]{partI}) to the theory of stacks over Segal sites presented here. Note that
Theorem \ref{tstrict} implies that
$$\hat{T}\simeq L(\mathsf{SSet}^{T})=L(\mathsf{SSet}^{T^{op}})=L(\mathsf{Pr}(T,\mathsf{SSet})).$$

The strictification theorem \ref{tstrict} has one fundamental consequence which is the \textit{Yoneda
lemma for Segal categories} proved in \cite{s3}. Let $A$ be a Segal category (say in
$\mathbb{U}$); we know that we can find
an $S$-category $A'$ and an equivalence $A' \simeq A$. For the $S$-category $A'$ one has
a natural simplicially enriched Yoneda morphism
$(A')^{op}\times A'  \longrightarrow \mathsf{SSet}_{\mathbb{U}}$. By adjunction this gives
a well defined morphism of $S$-categories
$A' \longrightarrow \mathsf{SSet}_{\mathbb{U}}^{(A')^{op}}$. Composing this with the morphism of Theorem
\ref{tstrict}, we get a morphism of Segal categories
$$h : A\simeq A' \longrightarrow \hat{A'}\simeq \hat{A},$$ well defined in
the homotopy category of Segal categories.
Now, Theorem \ref{tstrict}, results of \cite{dk4} and again the simplicially enriched Yoneda lemma together imply that
the morphism $h$ is fully faithful. This is the Yoneda lemma for the Segal category
$A$ of \cite{s3}. More generally, the same kind of argument proves the following proposition.

\begin{prop}\emph{(Yoneda lemma)}\label{Yoneda}
Let $A$ be a Segal category, and $a \in A$ an object.
Let $F \in \hat{A}$ be a morphism from $A^{op}$ to $\mathsf{Top}$. There exists a natural
equivalence of simplicial sets
$$F(a) \simeq \hat{A}_{(h_{a},F)}.$$
\end{prop}

\medskip

We finish this subsection with the following definition of representable and corepresentable functors.

\begin{df}\label{drep}
Let $A$ be a Segal category and $F \in \mathbb{R}\underline{Hom}(A^{op},\mathsf{Top})$
be an object. $F$ is called \emph{representable} if
it is equivalent to $h_{a}$ for some $a \in A$.

Dually, an object $F \in \mathbb{R}\underline{Hom}(A,\mathsf{Top})$ is called
\emph{corepresentable} if it is representable when considered as
a functor from $(A^{op})^{op}$ to $\mathsf{Top}$ (i.e. equivalent to
$h_{a}$ for some $a \in A^{op}$).
\end{df}

\end{subsection}

\begin{subsection}{Comparison}

Recall from \cite[\S 3.8]{partI} that a ($t$-complete) \textit{model topos}
is a model category Quillen equivalent to the model category $\mathsf{SPr}_{\tau}(T)$
of stacks over  some $S$-site $(T,\tau)$. The results of the previous subsection imply
that the notions of model topos and of Segal topos are related via the
functor $L$, which sends a model category to its simplicial localization
(along equivalences). More presicely, if one starts with
a $\mathbb{U}$-small $S$-site $(T,\tau)$, then one has the model
categories of pre-stacks $\mathsf{SPr}(T)$ and of stacks $\mathsf{SPr}_{\tau}(T)$
on $T$, as defined in \cite[\S 3]{partI}. One the other hand, the $S$-site $(T,\tau)$ might be
considered in a trivial way as a $\mathbb{U}$-small Segal site, and one can consider the
associated Segal categories of pre-stacks $\hat{T}$ and of stacks $T^{\sim,\tau}$.
Theorem \ref{tstrict} implies the following comparison result.

\begin{cor}\label{ccomp}
\begin{enumerate}
\item If $M$ is a $\mathbb{U}$-model topos then $LM$ is a $\mathbb{U}$-Segal topos.
\item
Let $(T,\tau)$ be a $\mathbb{U}$-small $S$-site. There exists natural isomorphisms in $\mathrm{Ho}(\mathsf{PrSeCat})$
$$L(\mathsf{SPr}(T))\simeq \hat{T} \qquad L(\mathsf{SPr}_{\tau}(T))\simeq T^{\sim,\tau}.$$
\end{enumerate}
\end{cor}

Therefore the theory of model topoi, investigated in \cite{partI}, embeds in the
present theory of Segal topoi where it actually finds a richer environment,
mainly because when viewed as Segal topoi (via the simplicial localization $L$),
two model topoi \textit{do} have a Segal category of geometric morphisms between
them (as explained in Definition \ref{dgeom}) while it is not clear what a model
category of geometric morphisms between model topoi should be.

Let us also mention that all results and constructions of \cite{partI}
provide through corollary \ref{ccomp} analogs results and constructions
in the context of Segal topoi. This way, we obtain
the existence of internal Hom's, a theory of truncated objects and truncation
functors \dots. \\

Finally, corollary \ref{ccomp} together with theorem \ref{cmp} implies the following important 
fact.

\begin{cor}\label{cpres}
A $t$-complete $\mathbb{U}$-Segal topos is
a presentable Segal category in the sense of Lurie (\cite[Prop. 1.4.1]{lu}).
\end{cor}

\textit{Proof:} Indeed, $T$ is equivalent to some $LSPr_{\tau}(T)$. As
$SPr_{\tau}(T)$ is a cofibrantly generated model category, $T$ is a presentable
Segal category (see \cite{s3} for details). \hfill $\Box$ \\

\begin{rmk}
\emph{Corollary \ref{cpres} states that any $t$-complete Segal topos is
a presentable Segal category. The question whether any Segal
topos is a presentable Segal category is open.}
\end{rmk}

\end{subsection}

\end{section}

\begin{section}{Further examples and applications}

In this last section we give some further materials and some examples of applications
of the theory of Segal topoi.

\begin{subsection}{Giraud's theorem}

In order to state the conjectural Giraud's theorem for Segal topoi, we need
some preliminaries notations and definitions. \\

\begin{itemize}

\item A simplicial set is called \emph{essentially $\mathbb{U}$-small} if it is
equivalent to a $\mathbb{U}$-small simplicial set.

\item A coproduct $\coprod_{i \in I}x_{i}$ in a Segal category $A$ with inital object
$\emptyset$ is called
\emph{disjoint} if for any $i\neq j$ the square
$$\xymatrix{\emptyset \ar[r] \ar[d] & x_{i} \ar[d] \\
 x_{j} \ar[r] & \coprod_{i\in I}x_{i}}$$
is cartesian in $A$.

\item Let $X_{*} : \Delta^{op} \longrightarrow A$ be a simplicial object
in a Segal category $A$. We say that $X_{*}$ is a \emph{groupoid object}
if it satisfies the following two conditions.
\begin{enumerate}

\item The Segal morphisms (defined as in \S 2)
$$X_{n} \longrightarrow \underbrace{X_{1}\times_{X_{0}}X_{1} \dots \times_{X_{0}}X_{1}}_{n \; times}$$
are equivalences (this condition includes the fact that the fibered product
of the left hand side exists in $A$).

\item The morphism
$$d_{2}\times d_{1} : X_{2} \longrightarrow X_{1}\times_{d_{0},X_{0},d_{1}}X_{1}$$
is an equivalence.

\end{enumerate}

\item Finally, let $A$ be a $\mathbb{V}$-small Segal category whose simplicial sets of morphisms
are essentially $\mathbb{U}$-small. A $\mathbb{U}$-small \textit{set of strong generators} for $A$ is
a $\mathbb{U}$-small set $E$ of objects in $A$ (considered as a full sub-Segal category of $A$), such that
the natural morphism
$$A \longrightarrow \hat{A} \longrightarrow \hat{E}$$
is fully faithful (the first morphism is the Yoneda embedding $h$ of Proposition \ref{Yoneda}, and
the second one is the restriction morphism from $A$ to $E$).

\end{itemize}

The analog of Giraud's theorem for Segal categories is the following.

\begin{conj}\emph{(Giraud's theorem for Segal topoi)}\label{giraud}
Let $A$ be a $\mathbb{V}$-small Segal category. Then $A$ is a $\mathbb{U}$-Segal topos if and
only if it satisfies the following conditions.

\begin{enumerate}
\item The simplicial sets of morphisms in $A$ are essentially $\mathbb{U}$-small.

\item The Segal category $A$ has all $\mathbb{U}$-small colimits, and coproducts in $A$ are
disjoint.

\item For all groupoid objects $X_{*}$ in $A$, the natural
morphism
$$X_{*} \longrightarrow N(X_{0} \rightarrow |X_{*}|)$$
is an equivalence of simplicial objects in $A$ (here $|X_{*}|$ is the colimit
of the simplicial diagram $X_{*}$ and $N(X_{0} \rightarrow |X_{*}|)$ is the nerve
of the natural morphism $X_{0} \rightarrow |X_{*}|$).

\item Colimits in $A$ are stable by pullbacks. In other words, for any $\mathbb{U}$-small
category $I$, any $I$-diagram $x_{*}$ in $A$, and any morphisms
$\xymatrix{\mathrm{Colim}_{i\in I}x_{i} \ar[r] & z & \ar[l] y}$, the natural morphism
$$\mathrm{Colim}_{i\in I}\left(
x_{i}\times_{z}y\right) \longrightarrow \left(\mathrm{Colim}_{i\in I}x_{i}\right)\times_{z}y$$
is an equivalence.

\item The Segal category $A$ has a $\mathbb{U}$-small set of strong generators.

\end{enumerate}
\end{conj}

Conjecture \ref{giraud} can be directly compared with Theorem $1$ of the appendix of \cite{mm}.
It can also be compared with the version of Giraud's theorem in the context of model topoi stated by C. Rezk
(unpublished). A proof of the conjecture \ref{giraud} have been recently given by J. Lurie in \cite{lu}
(at least for presentable Segal categories). \\

Finally, one can easily prove the following proposition, which together with
the previous conjecture would give a Giraud's theorem of $t$-complete Segal topoi, and
therefore, by Theorem \ref{cmp}, a characterization of Segal categories of the form $T^{\sim,\tau}$.

For this, let us recall that a morphism $f : F \longrightarrow G$
in a Segal topos $A$ is an \emph{epimorphism} if the natural morphism
$$|N(F\rightarrow G)| \longrightarrow G,$$
from the geometric realization of the nerve of $f$ down to $G$, is an
equivalence in $A$. We will then say that an augmented simplicial object $X_{*} \rightarrow Y$ in a Segal topos $A$
is a \emph{hyper-cover} if for any $n\geq 0$ the induced morphism
$$X_{*}^{\mathbb{R} \Delta^{n}} \longrightarrow
X_{*}^{\mathbb{R} \partial \Delta^{n}}\times_{Y^{\mathbb{R} \partial \Delta^{n}}}Y^{\mathbb{R} \Delta^{n}}$$
is an epimorphism.

\begin{prop}
A Segal topos $A$ is $t$-complete if and only if for all hyper-cover
$X_{*} \longrightarrow Y$ the induced morphism
$$|X_{*}| \longrightarrow Y$$
is an equivalence.
\end{prop}

\end{subsection}

\begin{subsection}{A Galois interpretation of homotopy theory}

In this number, we will give an application of the notion of geometric morphisms between Segal topoi to the
\textit{Galois interpretation of homotopy theory} of \cite{to3}. \\

Let $X$ be a CW complex in $\mathbb{U}$, and $\mathsf{SPr}(X)$ be
the category of $\mathbb{U}$-simplicial presheaves on $X$
(i.e. presheaves of $\mathbb{U}$-small simplicial sets on the site
of open subsets of $X$). The category $\mathsf{SPr}(X)$ is
endowed with its local projective model structure of \cite{bl}.
We consider the full subcategory $\mathsf{PrLoc}(X)$ of
$\mathsf{SPr}(X)$, consisting of \textit{locally $h$-constant objects} in the
sense of \cite{to3}. They are the simplicial presheaves
$F$ which are, locally on $X$, equivalent (for the local model structure) to a constant simplicial presheaf.
We define $\mathsf{Loc}(X):=L\mathsf{PrLoc}(X)$, the simplicial
localization of $\mathsf{PrLoc}(X)$ along the local equivalences, and
we view it as a $\mathbb{V}$-small Segal category.

The following theorem is the higher analog of the Galois interpretation of the theory of fundamental groupoids
explored by Grothendieck in \cite{sga1}. It is also an extension of the main theorem of \cite{to3} to the
case of unbased spaces.

\begin{thm}\label{tgal}
\begin{enumerate}
\item The Segal category $\mathsf{Loc}(X)$ is a $t$-complete $\mathbb{U}$-Segal topos.

\item Let $X$ and $Y$ be two $\mathbb{U}$-small $CW$-complexes. The Segal category
$\mathbb{R}\underline{Hom}^{geom}(\mathsf{Loc}(X),\mathsf{Loc}(Y))$ is a Segal groupoid (i.e.
its homotopy category is a groupoid, or equivalentely all its morphisms are equivalences).

\item Let $X$ and $Y$ be two $\mathbb{U}$-small $CW$-complexes. There exists a natural
isomorphism in $\mathrm{Ho}(\mathsf{SSet})$
$$\mathbb{R}\underline{Hom}(X,Y)\simeq |\mathbb{R}\underline{Hom}^{geom}(\mathsf{Loc}(X),\mathsf{Loc}(Y))|,$$
where $|\mathbb{R}\underline{Hom}^{geom}(\mathsf{Loc}(X),\mathsf{Loc}(Y))|$ is the geometric realization
(see \cite[\S 2]{sh})
of the Segal category $\mathbb{R}\underline{Hom}^{geom}(\mathsf{Loc}(X),\mathsf{Loc}(Y))$.

\item The functor
$$\begin{array}{cccc}
\mathsf{Loc} : & \mathrm{Ho}(\mathsf{Top}) & \longrightarrow
& \mathrm{Ho}(\mathsf{SeT}_{\mathbb{U}})\subset \mathrm{Ho}(\mathsf{PrSeCat}) \\
& X & \mapsto & \mathsf{Loc}(X)
\end{array}$$
is fully faithful.

\end{enumerate}
\end{thm}

\textit{Sketch of proof:} The theorem is mainly a consequence of results in \cite{to3}. \\

$(1)$ This follows from \cite[Theorem 2.13, 2.22]{to3} and the strictification theorem \ref{tstrict},
which shows that $\mathsf{Loc}(X)$ is equivalent to
$\widehat{BG}$, where $G$ is the simplicial group of loops on $X$
(we assume here that $X$ is connected for the sake of simplicity), and
$BG$ is the Segal category with a unique object and $G$ as the endomorphism simplicial monoid. \\

$(2)$ As above, let us write $\mathsf{Loc}(X)\simeq \widehat{BG}\simeq \mathbb{R}\underline{Hom}(BG^{op},\mathsf{Top})$.
By adjunction one has
$$\mathbb{R}\underline{Hom}^{geom}(\mathsf{Loc}(X),\mathsf{Loc}(Y))
\simeq \mathbb{R}\underline{Hom}(BG,\mathbb{R}\underline{Hom}^{geom}(\mathsf{Top},\mathsf{Loc}(Y)).$$
As $\mathsf{Loc}(*)\simeq \mathsf{Top}$,
this shows that one can assume $X=*$. One can also clearly assume that $Y$ is connected. Therefore, it is enough to
show that all objects in $\mathbb{R}\underline{Hom}^{geom}(\mathsf{Top},\mathsf{Loc}(Y))$ are equivalent
(compare with \cite[Corollary 5.7]{sga1}).

Let $y \in Y$ be a base point and $\omega_{y} : \mathsf{Loc}(Y) \longrightarrow \mathsf{Top}$
the fiber-at-$x$ morphism. We will prove that any
$\omega \in \mathbb{R}\underline{Hom}^{geom}(\mathsf{Top},\mathsf{Loc}(Y))$
is equivalent to $\omega_{y}$. Indeed, by \cite[Theorem 2.22]{to3} we know that
$\omega_{y}$ is co-represented by some $E \in \mathsf{Loc}(Y)$, corresponding to the action
of the loop group of $Y$ on itself. Since $E\neq \emptyset$, the exactness assumption on
$\omega$ implies that $\omega(E)\neq \emptyset$. The Yoneda lemma \ref{Yoneda} implies
that $\omega(E)\simeq \mathbb{R}\underline{Hom}^{geom}(\mathsf{Top},\mathsf{Loc}(Y))_{(\omega_{y},\omega)}$,
and therefore
one sees that there exists a morphism $u : \omega_{y} \rightarrow \omega$. We need to show that
$u$ is an equivalence, or equivalentely that for any $F \in \mathsf{Loc}(Y)$, the induced morphism
$u_{F} : F_{y} \rightarrow \omega(F)$ is an equivalence of simplicial sets. Considering this morphism on the
level of $\pi_{0}$, and using \cite[Corollary 6.3]{sga1}, one sees that
for any $F \in \mathsf{Loc}(Y)$ the induced morphism $\pi_{0}(u_{F}) : \pi_{0}(F_{y}) \rightarrow \pi_{0}(\omega(F))$
is bijective. Moreover, applying this to $F^{\mathbb{R} K}$, for $K$ a finite simplicial set,
and using exactness of $\omega$ we see that $F_{y} \rightarrow \omega(F)$ is in fact an equivalence. \\

$(3)$ As above, one can suppose that $X=*$ and that $Y$ is connected.
Let us write $Y \simeq BG$, where
$G$ is a simplicial group. Using \cite[Theorem 2.13.,2.22]{to3}, we need to show that the natural morphism
$$BG \longrightarrow \mathbb{R}\underline{Hom}^{geom}(\mathsf{Top},\widehat{BG})$$
is an equivalence of Segal categories. The fact that this morphism is fully faithful is
\cite[Corollary 3.2]{to3}. On the other hand, point $(2)$ above shows that all
objects in $\mathbb{R}\underline{Hom}^{geom}(\mathsf{Top},\widehat{BG})$
are equivalent, i.e. that the above morphism is also essentially surjective. \\

$(4)$ The functor
$$\mathsf{Loc} : \mathrm{Ho}(\mathsf{Top}) \longrightarrow \mathrm{Ho}(\mathsf{SeT}_{\mathbb{U}})$$
has a right adjoint, sending a Segal topos $A$ to the geometric realization of the Segal
groupoid $\mathbb{R}\underline{Hom}^{geom}(\mathsf{Top},A)^{int}$
(the \textit{interior} of the Segal category $\mathbb{R}\underline{Hom}^{geom}(\mathsf{Top},A)$, see
\cite[\S 2]{sh}. It is the maximal sub-Segal groupoid
of $\mathbb{R}\underline{Hom}^{geom}(\mathsf{Top},A)$). The fact that the adjunction
morphism $X \longrightarrow |\mathbb{R}\underline{Hom}^{geom}(\mathsf{Top},\mathsf{Loc}(X))|$
is an isomorphism for any $X$ in $\mathrm{Ho}(\mathsf{Top})$ is the content of point $(3)$. \hfill $\Box$ \\

>From Theorem \ref{tgal} $(4)$, one deduces the following reconstruction
result. Note that the Segal category $\mathsf{Top}=L\mathsf{SSet}_{\mathbb{U}}$ of simplicial
sets is equivalent to $\mathsf{Loc}(*)$, $*$ denoting the one-point space.

\begin{cor}
The Segal topos $\mathsf{Loc}(X)$ determines $X$. More precisely, there exists a natural isomorphism in
$\mathrm{Ho}(\mathsf{SSet})$
$$X\simeq |\mathbb{R}\underline{Hom}^{geom}(\mathsf{Top},\mathsf{Loc}(X))|.$$
In particular, if $X$ and $Y$ are two $\mathbb{U}$-small $CW$-complexes such that the Segal categories
$\mathsf{Loc}(X)$ and $\mathsf{Loc}(Y)$ are equivalent, then $X$ and $Y$ are homotopy equivalent.
\end{cor}

\end{subsection}

\begin{subsection}{Homotopy types of Segal sites}

The reconstruction of a homotopy type from the Segal category of locally constant stacks on it
explained in the last subsection suggests that one could try to extend homotopy theory
to Segal sites, in the same way as the theory of fundamental group is extended to the more
general setting of \cite[Exp. V]{sga1}. For nice enough Grothendieck sites, a pro-homotopy type has been
constructed in \cite{am}. However, the approach of \cite{am} is very different from
the galoisian point of view originally adopted in \cite[Exp. V]{sga1}, and no relations have still been
made explicit between certain stacks on the site and certain \textit{stacks on its associated homotopy type},
as, for example, locally constant sheaves corresponds to continuous representations
of the fundamental group. Moreover, the construction of \cite{am} only gives a pro-object in the
homotopy category of spaces, whereas various
works on \'etale homotopy theory and pro-finite completions suggest that one should
actually expect an object in the homotopy category of pro-spaces instead.

In this last part we propose a very general approach to define homotopy types
of Segal sites, which follows the original galoisian point of view of \cite{sga1} and \cite{to3}.
Not all the details of the construction will be presented here. \\

Let $T$ be a $\mathbb{U}$-Segal topos. We assume that $T$ is \textit{presentable Segal category
in the sense of \cite{lu}} (e.g. $T$ is a $t$-complete Segal topos, see
Corollary \ref{cpres}).

We wish to define an associated morphism of Segal categories
$H_{T} : \mathsf{Top} \longrightarrow \mathsf{Top}$.\\ For a $\mathbb{U}$-small simplicial set $Y$, we consider
$\mathsf{Loc}(Y)$, the Segal category of locally constant stacks on its geometric realization
(by \cite[Theorem 2.13]{to3} $\mathsf{Loc}(Y)$ is also equivalent to the comma Segal category $\mathsf{Top}/Y$ of objects
over $Y$) and the Segal category of geometric morphisms
$\mathbb{R}\underline{Hom}^{geom}(T,\mathsf{Loc}(Y))$.

On the other hand we can look at the unique geometric morphism $p : T \longrightarrow \mathsf{Top}$, with
inverse image $p^{*} : \mathsf{Top} \longrightarrow T$, and define the cohomology
of $T$ with coefficients in a $\mathbb{U}$-small simplicial set $Y$ as
$$\mathbb{H}(T,Y):=T_{(*,p^{*}(Y))}.$$
Note that, as the simplicial sets of morphisms in $T$ are essentially
$\mathbb{U}$-small, so is $\mathbb{H}(T,Y)$. We have the following comparison result.

\begin{lem}\label{lhom}
The Segal category
$\mathbb{R}\underline{Hom}^{geom}(T,\mathsf{Loc}(Y))$ is a Segal groupoid
whose geometric realization is naturally equivalent to the cohomology space of
$T$ with coefficients in $Y$
$$|\mathbb{R}\underline{Hom}^{geom}(T,\mathsf{Loc}(Y))| \simeq \mathbb{H}(T,Y).$$
\end{lem}

\textit{Proof:} Let us first describe the natural morphism
$$\phi_{T,Y} : |\mathbb{R}\underline{Hom}^{geom}(T,\mathsf{Loc}(Y))| \longrightarrow \mathbb{H}(T,Y).$$
First of all, by \cite[Theorem 2.13]{to3}, one has a natural equivalence
between $Loc(Y)_{(*,\underline{Y})}$ and $Map(Y,Y)$, where
$\underline{Y}$ is the constant stack on $Y$ with fiber $Y$. Therefore, the identity
provides a natural element $c_{Y} \in Loc(Y)_{(*,\underline{Y})}$.
Now, if $f^{*} : Loc(Y) \longrightarrow T$ is the inverse image of a
geometric morphism $f \in \mathbb{R}\underline{Hom}^{geom}(T,\mathsf{Loc}(Y))$,
the image of $c_{Y}$ by $f^{*}$ gives an element in
$T_{(*,f^{*}(\underline{Y}))}\simeq \mathbb{H}(T,Y)$. This defines the morphism
$\phi_{T,Y}$.

To prove that $\mathbb{R}\underline{Hom}^{geom}(T,\mathsf{Loc}(Y))$ is a Segal groupoid
and that $\phi_{T,Y}$ is an equivalence, one writes $T$ as a left localization
of $\hat{B}$, for some Segal category $B$. Then,
$\mathbb{R}\underline{Hom}^{geom}(T,\mathsf{Loc}(Y))$ is a left exact localization
of $\mathbb{R}\underline{Hom}^{geom}(\hat{B},\mathsf{Loc}(Y))$, and
using the functoriality of $\phi_{T,Y}$ in $T$ one is reduced to the case
where $T$ is of the form $\hat{B}$ where the lemma follows easily from
the adjunction formula
$$\mathbb{R}\underline{Hom}^{geom}(\hat{B},\mathsf{Loc}(Y))\simeq
\mathbb{R}\underline{Hom}(B^{op},\mathbb{R}\underline{Hom}^{geom}(Top,\mathsf{Loc}(Y))),$$
and theorem \ref{tgal}. \hfill $\Box$. \\

We define $H_{T}(Y)$ as the geometric realization of the Segal groupoid
$\mathbb{R}\underline{Hom}^{geom}(T,\mathsf{Loc}(Y))$:
$$H_{T}(Y):=|\mathbb{R}\underline{Hom}^{geom}(T,\mathsf{Loc}(Y))| \in \mathsf{Top}.$$

This defines a morphism $H_{T} : \mathsf{Top} \longrightarrow \mathsf{Top}$, which is abstractly considered
as ``some'' homotopy type associated to the Segal topos $T$. By the previous lemma, one
also has a natural equivalence $H_{T}(Y)\simeq \mathbb{H}(T,Y)$, for any $\mathbb{U}$-small simplicial set $Y$.

\begin{df}\label{dhom}
The \emph{homotopy shape} of the Segal topos $T$ is defined to be the
object $H_{T} \in \mathbb{R}\underline{Hom}(\mathsf{Top},\mathsf{Top})$ defined above.
\end{df}

A fundamental consequence of Theorem \ref{tgal} is that if $T$ is the Segal topos of stacks over
the Grothendieck site of a CW complex $X$, then $H_{T}$ is corepresented by the homotopy type of $X$.
In general we do not
expect $H_{T}$ to be corepresentable, but one can prove that
$H_{T}$ is \textit{pro-corepresentable in the context of Segal categories}. This
last fact is more or less equivalent to the fact that the morphism
$H_{T}$ is left exact, which follows from lemma \ref{lhom}. Let us be more precise.

A Segal category $K$ is called \textit{finite}, if for any
filtered system of Segal categories $\{B_{i}\}_{i\in I}$, the
natural morphism
$$Colim_{i\in I}\mathbb{R}\underline{Hom}(K,B_{i})
\longrightarrow \mathbb{R}\underline{Hom}(K,Colim_{i\in I}B_{i})$$
is an equivalence. A Segal category $A$ is \textit{left filtered} if
for any finite Segal category $K$ and any $K$-diagram $x \in \mathbb{R}\underline{Hom}(K,A)$,
there exists an object $a \in A$ and a morphism
$a \rightarrow x$ in $\mathbb{R}\underline{Hom}(K,A)$ (here, $a$ also denotes
the constant $K$-diagram in $A$ with value $a$).

\begin{prop}\label{phom}
Let $T$ be a $\mathbb{U}$-Segal topos which
is a presentable Segal category (e.g. a $t$-complete Segal topos).
There exists a left filtered
$\mathbb{U}$-small Segal category $A$
and an $A$-diagram $\chi(T) \in \mathbb{R}\underline{Hom}(A,\mathsf{Top})$ such that
for any $Y\in Top$ there is a natural equivalence
$$H_{T}(Y)\simeq Colim_{a\in A^{op}}\underline{Hom}(\chi(T)_{a},Y).$$
In the other words, the two endomorphisms of $\mathsf{Top}$
$$Y \mapsto H_{T}(Y) \qquad Y \mapsto Colim_{a\in A^{op}}\underline{Hom}(\chi(T)_{a},Y)$$
are equivalent as objects in $\mathbb{R}\underline{Hom}(\mathsf{Top},\mathsf{Top})$.
\end{prop}

\textit{Sketch of proof:} Let us recall that we have
assumed $T$ to be a presentable Segal category. Therefore, there exists a
regular cardinal $\alpha \in \mathbb{U}$ such that, for any simplicial set $Y$,
$H_{T}(Y)\simeq Colim_{Y_{i}}H_{T}(Y_{i})$, where $Y_{i}$ runs
through the sub-simplicial sets of $Y$ whose cardinality is less than
$\alpha$ (chose $\alpha$ so that $*$ is a $\alpha$-small object in $T$ and
observe that $H_{T}(Y)\simeq T_{(*,\underline{Y})}$).

We now consider $H_{T}$, restricted to the full sub-category
$\mathsf{Top}_{\leq \alpha}$ of $\mathsf{Top}$ consisting of simplicial sets of
cardinality less than $\alpha$. We define $A$ to be the Segal
category of objects of the morphism $H_{T} : \mathsf{Top}_{\leq \alpha}
\longrightarrow \mathsf{Top}$. In other words, $A$ is the Segal category of pairs
$(Y,x)$, where $Y \in \mathsf{Top}_{\leq \alpha}$ and $x\in H_{T}(Y)$ (the
Segal category $A$ is also denoted by $\int H_{T}$ in \cite{sh}).
As $H_{T}$ is left exact and by the choice of the cardinal $\alpha$, one
sees that $A$ satisfies the conditions of the proposition. \hfill $\Box$ \\

The diagram $\chi(T)$ of proposition \ref{phom} is called the
\textit{pro-homotopy type of the Segal topos $T$}. It is not a pro-simplicial in
the usual sense, as the indices live in a filtered Segal category $A$ rather than
in an actual filtered category. However, this seems to be not of so much
importance, as pro-objects in the Segal setting behave very much the same
way as usual pro-objects. Also, as argued in \cite{lu}, it seems
that the two notions are more or less equivalent, and that the
Segal category $A$ of proposition \ref{phom} can in fact be chosen to be
a usual filtered category.

In any case, when $A$ is a filtered Segal category then its homotopy category
$\mathrm{Ho}(A)$ is a filtered category in the usual sense. Therefore, the $A$-diagram
$\chi(T)$ of proposition
\ref{phom} induces a $\mathrm{Ho}(A)$-diagram in $\mathrm{Ho}(\mathsf{Top})$, i.e. a pro-object in the homotopy category of spaces. When
$T$ is the Segal topos of stacks over a locally connected
Grothendieck site, we suspect that this induced pro-object is equivalent to the
one constructed by Artin and Mazur in \cite{am}. The pro-homotopy type $\chi(T)$
of the Segal topos $T$ is therefore a refinement and a generalization
of Artin-Mazur's construction. \\

An interesting example of application of homotopy types of Segal topoi would be the study of the
\textit{\'etale homotopy type of the sphere spectrum}, defined by using the
Segal site $(Spec\, \mathbb{S},\textrm{\'et})$ defined in \cite[\S 5]{partI}. We can, for example, ask the following more
precise question.

\begin{Q}\label{quest2}
Let $T:=Spec\, \mathbb{S}^{\sim,\textrm{\'et}}$ be the small \'etale Segal topos of the sphere spectrum
defined in \cite[\S 5]{partI}, and $\mathbb{S}^{\wedge,\textrm{\'et}}:=\chi(T)$ be its pro-homotopy type
as defined above. Describe $\mathbb{S}^{\wedge,\textrm{\'et}}$, and in particular compare it
with the etale homotopy type of $Spec\, \mathbb{Z}$.
\end{Q}

One could also ask, for each rational prime $p$, a similar
question for the $p$-localized sphere spectrum $\mathbb{S}_{p}$;
in this case it seems natural to ask whether
the corresponding pro-corepresentative space is determined
by the Morava $K$-theories $K(n)$ or $E$-theories $E(n)$, $n\geq 0$, in analogy with the standard chromatic picture.\\

\end{subsection}

\end{section}


\begin{thebibliography}{90}


\bibitem[SGA1]{sga1} M. Artin, A. Grothendieck,
\textit{Rev\^etements \'etales et groupe fondamental}, Lecture notes in Math.
$\mathbf{224}$, Springer-Verlag 1971.

\bibitem[SGA4-I]{sga4} M. Artin, A. Grothendieck, J. L. Verdier, \textit{Th\'eorie des topos et cohomologie
\'etale des sch\'emas- Tome 1}, Lecture Notes in Math \textbf{269}, Springer Verlag, Berlin, 1972.

\bibitem[A-M]{am} M. Artin, B. Mazur, \textit{\'{E}tale homotopy},
Lecture Notes in Math. $\mathbf{100}$, Springer-Verlag $1969$.

\bibitem[Bl]{bl} B. Blander, \textit{Local projective model structure on simpicial presheaves},
$K$-theory \textbf{24} (2001) No. $3$, 283-301.

\bibitem[Cis]{cis} D-C. Cisinski,
\textit{Th\'eories homotopiques dans les topos}, to appear in JPAA.

\bibitem[Du]{du2} D. Dugger \textit{Combinatorial model categories have presentations},
Adv. in Math. $\mathbf{164}$ (2001), 177-201.

\bibitem[D-K$1$]{dk1} W. Dwyer,  D. Kan, \textit{Simplicial localization of categories},
J. Pure and Appl. Algebra $\mathbf{17}$ (1980), 267-284.

\bibitem[D-K$2$]{dk3} W. Dwyer, D. Kan, \textit{Homotopy commutative diagrams and their realizations}, J. Pure Appl. Algebra
\textbf{57} (1989) No. 1, 5-24.

\bibitem[D-K$3$]{dk4} W. Dwyer,  D. Kan, \textit{Function complexes in homotopical algebra},
Topology $\mathbf{19}$ $(1980)$, $427-440$.

\bibitem[Fr]{fr} E. Friedlander, \textit{\'Etale homotopy type of simplicial schemes},
Annals of Mathematics Studies, $\mathbf{104}$. Princeton University Press,
Princeton, N.J., 1982.

\bibitem[Gr]{gr} A. Grothendieck, \textit{Pursuing stacks}, unpublished manuscript.

\bibitem[Hi]{hi} P. Hirschhorn, \textit{Localization of model categories}, Book in preparation available
at \textsf{http://www-math.mit.edu/$^{\sim}$psh}

\bibitem[H-S]{sh} A. Hirschowitz, C. Simpson, \textit{Descente pour les $n$-champs},
preprint available at \textsf{math.AG/$9807049$}.

\bibitem[Ho]{ho} M. Hovey, \textit{Model categories}, Mathematical surveys and monographs, Vol. $\mathbf{63}$,
Amer. Math. Soc., Providence 1998.

\bibitem[Ja]{ja} J. F. Jardine, \textit{Simplicial presheaves}, J. Pure and Appl. Algebra $\mathbf{47}$ (1987),
35-87.

\bibitem[Jo]{jo2} A. Joyal, \textit{Quasi-categories}, unpublished manuscript.

\bibitem[Le]{le} T. Leinster, \textit{A survey of definitions of $n$-category},
Theory and Applications of Categories, \textbf{10} (2002), 1-70.

\bibitem[Lu]{lu} J. Lurie, \textit{On $\infty$-topos}, preprint math.CT/0306109.

\bibitem[M-M]{mm} S. Mac Lane, I. Moerdijk, \textit{Sheaves in Geometry and Logic}, Springer, New York 1992.

\bibitem[P]{p} R. Pellissier, \textit{Cat\'egories enrichies faibles}, Th\`ese,  Universit\'e de Nice-Sophia Antipolis,
June 2002, available at \textsf{http://math.unice.fr/$^{\sim}$lemaire}.

\bibitem[Re]{rez2} C. Rezk, \textit{The notion of a homotopy topos}, unpublished, February $2001$.

\bibitem[Sc-Vo]{sv} R. Schw\"anzl, R. Vogt, \textit{Homotopy homomorphisms and the hammock localization},
Bol. Soc. Mat. Mexicana (2)
\textbf{37} (1992), No. 1-2, 431-448.

\bibitem[Si]{s3} C. Simpson, \textit{A Giraud-type characterisation of the simplicial categories
associated to closed model categories as $\infty$-pretopoi}, Preprint \textsf{math.AT/$9903167$}.

\bibitem[To]{to3} B. To\"en, \textit{Vers une interpr\'etation Galoisienne
de la th\'eorie de l'homotopie}, Cahiers de Top. et G\'eom. Diff. Cat. 43 (2002),
No. 4, 257-312.

\bibitem[To-Ve$1$]{web} B. To\"en, G. Vezzosi, \textit{Algebraic geometry over model categories.
A general approach to Derived Algebraic Geometry},
preprint available at \textsf{math.AG/$0110109$}.

\bibitem[To-Ve$2$]{partI} B. To\"en, G. Vezzosi, \textit{Homotopical algebraic geometry I: Topos theory},
Preprint \textsf{math.AG/$0207028$}; submitted.


\end{thebibliography}
\end{document}